\documentclass[11pt]{amsart}

\usepackage[foot]{amsaddr}
\usepackage{amsfonts,epsfig,fancyhdr,stmaryrd}
\usepackage{newlfont,amsfonts,amssymb,amsmath,amsthm,amsgen,amscd,datetime,dsfont,hyperref,tensor}
\usepackage[small,nohug,heads=littlevee]{diagrams}
\usepackage[normalem]{ulem}
\diagramstyle[labelstyle=\scriptstyle]
\usepackage{multicol,picinpar,enumerate,cleveref,mathabx}
\usepackage{euscript,color}

\newcommand{\tM}{\widetilde{M}}
\newcommand{\tg}{\widetilde{g}}

\newcommand{\Hol}{\mathrm{Hol}}
\newcommand{\+}{\oplus}

\renewcommand{\O}{\mathbf{O}}

\renewcommand{\span}{\mathrm{span}}
\newcommand{\hook}{\makebox[7pt]{\rule{6pt}{.3pt}\rule{.3pt}{5pt}}\,}

\newcommand{\1}{\mathbf{1}}
\newcommand{\e}{\mathrm{e}}
\renewcommand{\d}{\mathrm{d}}

\newcommand{\R}{\mathbb{R}} 
\newcommand{\N}{\mathbb{N}} 

\theoremstyle{definition}
\newtheorem{definition}{Definition}

\theoremstyle{plain}

\newtheorem*{lem*}{Lemma}

\newtheorem{theorem}[definition]{Theorem}
\newtheorem*{theorem*}{Theorem}

\numberwithin{equation}{section}

\begin{document}

\title
{Completeness of certain compact Lorentzian locally symmetric  spaces}

\thanks{This work was supported by 
 the Australian Research
Council (Discovery Program DP190102360).
 }
\author{Thomas Leistner}\address[Thomas Leistner, corresponding author]{School of Mathematical Sciences, University of Adelaide, SA~5005, Australia}\email{thomas.leistner@adelaide.edu.au}
\author{Thomas Munn}\address[Thomas Munn]{Lund University, 
    Faculty of Science,
    Centre for Mathematical Sciences,
    Box 118,
    22100~Lund,
    Sweden
}
\email{thomas.munn@math.lu.se}


\begin{abstract}
We show that a compact Lorentzian locally symmetric space is geodesically complete if the Lorentzian factor in the  local de\,Rham-Wu decomposition is  of Cahen-Wallach type or if the maximal flat factor is one-dimensional and time-like. Our proof uses a recent result by Mehidi and Zeghib and an earlier result by Romero and S\'{a}nchez. 
 \end{abstract}

\maketitle
\setcounter{tocdepth}{1}

\section{Introduction and statement of result}
In contrast to Riemannian manifolds, a compact Lorentzian manifold is not necessarily geodesically complete --- the fundamental example is the incomplete Clifton-Pohl torus, see \cite[Example 7.16]{oneill83}. Hence,  the question under which conditions completeness holds   is fundamental in Lorentzian geometry, especially when addressing classification problems for compact manifolds.

One condition that ensures that a compact Lorentzian manifold $(M,g)$ is complete is  a sufficient amount of global symmetry: Marsden showed that any compact homogeneous semi-Riemannian manifold is complete  \cite{marsden72}. This includes  symmetric spaces, i.e.~manifolds for which each point $p\in M$ is a fixed point of an isometry whose differential at $p$ is $-\mathrm{Id}$.  In a striking generalisation of Marsden's result in the Lorentzian context, Romero and S\'{a}nchez showed that $(M,g)$ is complete if it admits a   time-like conformal Killing vector field \cite{romero-sanchez95}, so a rather small amount of global symmetry.
 In the case of only local homogeneity, the question is much harder to answer,  and only when the dimension of $M$ is $3$, Zeghib and Dumitrescu \cite{zeghib-dumitrscu10} were able to show that $(M,g)$ is complete if it is locally homogeneous.
 
Other authors have imposed conditions on the curvature of compact Lorentzian manifold $(M,g)$:  if $(M,g)$ is flat \cite{carriere89} or, more generally,   has  constant sectional curvature \cite{klingler96}, then it is complete. In a combination of global and curvature assumptions, in \cite{leistner-schliebner13} it was shown that $(M,g)$ is complete if it is a pp-wave, i.e.~if $(M,g)$ admits a parallel null vector field and satisfies a certain curvature condition. Very recently this result was generalised by Mehidi and Zeghib  \cite{MehidiZeghib22}, who showed that the curvature condition can be dropped, i.e.~that the existence of a parallel null vector field on a compact Lorentzian manifold implies that the manifold is complete. These results  show that compactness and certain types of Lorentzian special holonomy imply completeness. Another result that  links completeness to the holonomy group was given in  \cite{Ake-HauSanchez16}. 

Of course, assumptions on the homogeneity  and on the curvature  are not independent of  each other. Examples are {\em locally symmetric spaces}, which are defined by the condition that the curvature tensor is parallel. A locally symmetric space  $(M,g)$ is locally isometric to a specific  symmetric space, the model space for $(M,g)$. Since  symmetric spaces are complete, by the de\,Rham-Wu decomposition theorem \cite{derham52,wu64}, their universal cover decomposes into indecomposable symmetric spaces, i.e.~those that cannot be further decomposed into a semi-Riemannian product. Hence, every locally symmetric space has a local product decomposition into  locally symmetric spaces that  are locally isometric to an indecomposable symmetric space. We refer to this as the {\em local de\,Rham-Wu decomposition} of a locally symmetric space, see Section \ref{sec-prod} for more details.
It has been shown  by Cahen and Wallach \cite{cahen-wallach70} that there are only two types of indecomposable Lorentzian symmetric spaces: the spaces of  nonzero constant sectional curvature ---  de Sitter and anti-de Sitter space --- and the  Cahen-Wallach spaces, see Section~\ref{sec-prod} for their definition.

If one could show that all compact locally symmetric spaces are complete, then they would be isometric quotients of their model spaces and  their classification could be reduced to the classification of compact Clifford-Klein forms, and hence to the algebraic problem of classifying co-compact isometric group actions (see for example the  results on compact quotients of Cahen-Wallach spaces and further references on  compact Lorentzian Clifford-Klein forms in \cite{KathOlbrich15}).
However, so far this has only been achieved in dimension $3$ (in \cite{zeghib-dumitrscu10} as a special case of locally homogeneous spaces) and for the indecomposable spaces: 
 the results in \cite{carriere89}, \cite{klingler96} and \cite{leistner-schliebner13} imply that a compact, indecomposable Lorentzian locally symmetric   space is geodesically complete. In the second author's  thesis \cite{tom-thesis} this was generalised to locally symmetric spaces which are local products of  Euclidean space and a Cahen-Wallach space. Using the same method as in \cite{tom-thesis},  the new result in \cite{MehidiZeghib22}, as well as the result in \cite{romero-sanchez95}, here we generalise this further.

\begin{theorem}
\label{thm}
Let $(M,g)$ a compact Lorentzian locally symmetric space whose local de\,Rham-Wu decomposition satisfies one of the following conditions:
\begin{enumerate}
\item the Lorentzian factor is  of Cahen-Wallach type, or
\item the maximal flat factor is one-dimensional and time-like.
\end{enumerate}
Then the time-orientable cover of $(M,g)$ admits a parallel  vector field and hence $(M,g)$ is geodesically complete.
\end{theorem}
This leaves open  the cases when $(M,g)$ does not have constant sectional curvature but 
 the  Lorentzian factor
 in
the local de\,Rham-Wu decomposition of $(M,g)$ does. The simplest examples of this situation are those when $(M,g)$ is locally isometric to a product of Minkowski space and a sphere, or to a product of (anti-)de Sitter space with Euclidean space or a sphere.
 In \cite{tom-thesis} an attempt was made to apply the approach in 
\cite{carriere89} and \cite{klingler96} to such cases, but the methods are too specific to the constant curvature model spaces and do not easily transfer to products thereof. 

\subsubsection*{Acknowledgements.}
A special case of Theorem \ref{thm} was obtained in the second author's MPhil thesis \cite{tom-thesis}, written  under supervision of the first author. 
We would like to thank Michael Eastwood for taking the role of co-supervisor and for helpful discussions and comments.

\section{The local de\,Rham-Wu decomposition of a Lorentzian locally symmetric space}\label{sec-prod}
In this section we will clarify what we mean by the {\em local de\,Rham-Wu decomposition} of a locally symmetric space and we will recall the notion of a Cahen-Wallach space.

We call a symmetric space $(M,g)$  {\em indecomposable} if the connected component of its holonomy group does not admit a non-trivial invariant subspace of the tangent space on which the metric remains non-degenerate. Since  symmetric spaces are complete, by the de\,Rham-Wu decomposition theorem \cite{derham52,wu64}, this is equivalent to the condition that the universal cover $(\tM,\tg)$ of $(M,g)$ does not decompose into a semi-Riemannian product. 
Hence, the universal cover of a symmetric space is globally isometric to a semi-Riemannian product of simply connected indecomposable symmetric spaces. 

Cahen and Wallach \cite{cahen-wallach70} have shown that an {\em indecomposable} Lorentzian symmetric space $(M,g)$ either has nonzero constant sectional curvature or is universally covered by $\tM=\R^{n+2} $ with 
\begin{equation}\label{CWmetric}
\tg= 2\d t \d v + x^a Q_{ab} x^b \d t^2 + \delta_{ab} \d x^a \d x^b,\end{equation}
where $Q_{ab}$ is a symmetric  $n\times n$-matrix with nonzero determinant and 
$(t,v,x^1,\ldots, x^n)$ are global coordinates on $\tM=\R^{n+2}$. Then $(\tM,\tg)$ is  called a {\em Cahen-Wallach space}. A Cahen-Wallach space admits a parallel null vector field $\frac{\partial}{\partial v}$, and the only non-vanishing terms (up to symmetry) of its curvature tensor are
\begin{equation}\label{CWcurv} R(\tfrac{\partial}{\partial t},\tfrac{\partial}{\partial x^a},\tfrac{\partial}{\partial t},\tfrac{\partial}{\partial x^b})=Q_{ab},\end{equation}
 for $a,b=1,\ldots, n$. The Ricci tensor $\rho$, as an endomorphism of $T\tM$, is
\begin{equation}\label{CWric} \rho= q\, \tfrac{\partial}{\partial v} \otimes \d t,\end{equation}
where $q$ is the trace of the matrix $Q_{ij}$. In particular, $\rho^2=0$. Two Cahen-Wallach spaces are isometric if their $Q_{ij}$'s have the same spectrum with the same multiplicities.

Applying the de\,Rham-Wu decomposition to the universal cover of a Lorentzian symmetric space $(M,g)$, we obtain that $(\tM,\tg)$ is globally isometric to a product of 
\begin{itemize}
\item[(a)]
an indecomposable Lorentzian symmetric space, i.e. of nonzero constant sectional curvature or a Cahen-Wallach space, or $(\R,-\d t^2) $,
\item[(b)]  irreducible Riemannian symmetric spaces,  and
\item[(c)]  Euclidean space.
\end{itemize}
Of course,  (b) or (c)  may not be present in the decomposition.  Moreover, any such a decomposition has a maximal flat factor, which is either empty or equal to  $(\R,-\d t^2) $, Minkowski space, or Euclidean space.

A locally symmetric space $(M,g)$  is a semi-Riemannian manifold with parallel curvature tensor, $\nabla R=0$. Every locally symmetric space is locally isometric to a globally symmetric space, and hence locally isometric to a product of indecomposable symmetric spaces. We refer to this local product decomposition as the {\em local de\,Rham-Wu decomposition} of $(M,g)$. We call a locally symmetric space {\em indecomposable} if it is locally isometric to an indecomposable symmetric space. Consequently, an indecomposable Lorentzian locally symmetric space either has constant sectional curvature or is locally isometric to a Cahen-Wallach space. 
The latter means that there are local coordinates such that the metric is of the form (\ref{CWmetric}).
Indecomposable compact Lorentzian locally symmetric spaces have been shown to be complete, the ones with constant sectional curvature in 
\cite{carriere89,klingler96} and the ones of Cahen-Wallach type in \cite{leistner-schliebner13}.

Since a locally symmetric Lorentzian symmetric space $(M,g)$ is locally isometric to a product of of symmetric spaces as in (a), (b), and (c), the tangent space at a point decomposes into orthogonal subspaces that correspond to the tangent spaces of the local factors. Moreover,  the curvature tensor  and the Ricci tensor of $(M,g)$ decompose accordingly  into the sum of the curvature tensors and Ricci tensors of the factors in the local decomposition.
Since irreducible Riemannian symmetric spaces of dimension $>1$ are Einstein with nonzero Einstein constant, this implies the following: if the local de\,Rham-Wu decomposition of $(M,g)$ does not contain a Cahen-Wallach factor, then the Ricci tensor of $(M,g)$ is given by the sum of multiples of the identity, and if it does contain a Cahen-Wallach factor, then the Ricci tensor is a sum of $\rho$ as in (\ref{CWric}) and multiples of the identity.

\section{Proof of Theorem \ref{thm}} 
Let $(M,g)$ be a locally symmetric Lorentzian space of dimension $m$ and let $\Hol_p$ be its holonomy group at a point $p\in M$.  Let $R$ be the curvature tensor of $(M,g)$ and $\rho$ the Ricci tensor considered as an endomorphism of the tangent bundle. Since $\nabla R=0$,  $\nabla \rho=0$, so both $R$ and $\rho$ are invariant under $\Hol_p$. 
We define two vector distributions, the {\em nullity of $R$},
\[ \mathcal N:=\{X\in TM\mid X\hook R=0\},\]
and the {\em generalised kernel of $\rho$},
\[ \mathcal K:=\{X\in TM\mid \exists m\in \N: \rho^m(X)=0\}.\]
Clearly, both $\mathcal N$ and $\mathcal K$ are invariant under $\Hol_p$, since for $H\in \Hol_p$ we have \[R(H \mathcal N|_p, X,Y,Z)=R(\mathcal N|_p, H^{-1}X, H^{-1}Y,H^{-1}Z)=0\] 
and \[
\rho(HX)=H^{-1}\rho(X),\]
for all $X,Y,Z\in T_pM$. This implies that also the orthogonal distributions $\mathcal N^\perp$ and $\mathcal K^\perp$ are invariant under $\Hol_p$. We clearly have that $\mathcal N\subset \mathcal K$. 

The remarks at the end of the previous section imply that in the case when the local de\,Rham-Wu decomposition of $(M,g)$ {\em does not} contain a Cahen-Wallach factor, we  have that $\mathcal N=\mathcal K$, $TM=\mathcal N\+\mathcal N^\perp$, and  $\cal N\not=0$ only if the local de\,Rham-Wu decomposition contains a flat factor. If the flat factor is one-dimensional and time-like, we have that $\mathcal N$ is a time-like line bundle and  an arbitrary element $H\in \Hol_p$ in the decomposition $T_pM=\mathcal N|_p\+\mathcal N^\perp|_p$ is of the form
\[H=\begin{pmatrix}\pm1 &0\\ 0& A\end{pmatrix},
\]
with $A\in \O(m-1)$.
Hence, the time-orientable cover of $(M,g)$ admits a time-like parallel vector field and is compact if $M$ is compact. A parallel vector field is of course conformally Killing, so we can apply the result in \cite{romero-sanchez95} to conclude that $(M,g)$ is complete. This proves case (2) in Theorem \ref{thm}.

From now on assume that the local de\,Rham-Wu decomposition of $(M,g)$ does contain a Cahen-Wallach factor. 
Since $\tfrac{\partial}{\partial v}$ is a local parallel null vector field,
 and with (\ref{CWcurv}), we get that
\[\mathcal L\ :=\  \mathcal N\cap \mathcal N^\perp\ \subset\ TM\]
is not trivial, and hence a bundle of null lines that is invariant under  the holonomy group of $(M,g)$. 
In particular, we have the following holonomy invariant filtrations
\[\mathcal L\ \subset\  \mathcal N\ \subset\  \mathcal K,\qquad \mathcal L\ \subset\  \mathcal N^\perp\cap \mathcal K \ \subset\  \mathcal K.\]
Moreover, since the Ricci-tensor of a Cahen-Wallach space is nilpotent and the irreducible Riemannian factors are Einstein with non-zero Einstein constant, $\mathcal K^\perp$  is non-degenerate, so that
\[TM=\mathcal K\+\mathcal K^\perp\]
is a holonomy invariant decomposition of $TM$.
At $p\in M$, we have  the $\Hol_p$-invariant decomposition
\[T_pM= \mathcal K|_p \+\mathcal K^\perp|_p ,\]
with  $\mathcal K|_p=T\+\overline{T}$, where $T\not=\{0\}$ is the
tangent space at $p$ of the local Cahen-Wallach factor and $\overline{T}$ the (possibly trivial) tangent space of the Euclidean (i.e. flat and Riemannian) factor. Note that the decomposition $\mathcal K|_p=T\+\overline{T}$  {\em may not} be $\Hol_p$-invariant. 
In order to describe the action of $\Hol_p$ on $T_pM$, we denote by $(t,v,x^1,\ldots x^n)$ the coordinates on the local Cahen-Wallach factor as in (\ref{CWmetric}) and fix the following basis of $\mathcal K|_p$,
\[e_-=\tfrac{\partial}{\partial v}, \quad e_1,\  \ldots ,\  e_n,\quad e_{n+1},\  \ldots,\  e_{n+k}, \quad e_+=\tfrac{\partial}{\partial v}-x^a Q_{ab} x^b\tfrac{\partial}{\partial t},\]
where $e_a=\tfrac{\partial}{\partial x^a}$ for $a=1,\ldots, n$   and $e_{n+b}$ for $b=1,\ldots, k$ an orthonormal basis for the Euclidean metric on $\overline{T}$.  In such a basis, the metric restricted to $\mathcal K|_p$  is of the form 
\[g|_p=\begin{pmatrix}0&0&1\\0&\1& 0 \\1&0&0\end{pmatrix}.\] 
In particular we have  
\[T=\span\{e_-,e_1\ldots, e_n,\e_+\}, \qquad \overline{T}=\span\{e_{n+1}\ldots, e_{n+k}\},\]
and that 
\[
\mathcal L|_p=\R\cdot e_-, \quad
(\mathcal N^\perp\cap \mathcal K) |_p=\span \{e_-,e_{1},\ldots, e_{n}\},
\quad
\mathcal N|_p=\span\{e_-,e_{n+1},\ldots, e_{n+k}\},\] 
are $\Hol_p$-invariant.
Then, in the $\Hol_p$-invariant decomposition $T_pM= \mathcal K|_p \+\mathcal K^\perp|_p$,
an arbitrary element $H\in \Hol_p$ is of the form
\[H=\begin{pmatrix}F&0\\0&G\end{pmatrix},\]
 where $F$, in the above basis, writes as
 \[F=\begin{pmatrix}\lambda&\ast&\ast
 &\ast
 \\0&A&0& u \\0&0&B&v\\
 0&0&0&\lambda^{-1}\end{pmatrix},\] 
with $A\in \O(n)$, $B\in \O(k)$, $u\in \R^{n}$, $v\in \R^{k}$,   $0\not= \lambda\in\R$, and the  $\ast$-terms are determined by the other terms.
In particular, we have that 
\[
H(e_a)=F(e_a)=A_a{}^be_b
 \mod \e_-, \qquad H(e_+)=F(e_+)=\lambda e_+ + u^ae_a +v^{\overline{c}}e_{\overline{c}}\mod e_-,\]
where  the unbarred indices $a,b,c\ldots$  run from $1$ to $n$ and the barred indices 
$\overline{a},\overline{b},\overline{c}, \ldots$
from $n+1$ to $n+k$.
Then, with 
(\ref{CWcurv}) and 
since $e_-\hook R=0$ and $e_{\overline c}\hook R=0$, we have
\[
(H^{-1}\cdot R)(e_a,e_+,e_b,e_+)
=
R(H e_a,H e_+,He_b,He_+)
=
R(Fe_a,F e_+,Fe_b,Fe_+)
=
\frac{A_a{}^cQ_{cd} A_b{}^d}{\lambda^2}.
\]
Hence, the holonomy invariance of $R$,  \[(H^{-1}\cdot R)(e_a,e_+,e_b,e_+)= R(e_a,e_+,e_b,e_+)=Q_{ab},\] yields the matrix equation for $Q=(Q_{ab})$, 
\[ \lambda^{2}Q=A^\top Q A.\]
Since $A\in \O(n)$ and $\det(Q)\not=0$, this implies that $\lambda=\pm 1$.
This however shows that the time-orientable cover of $(M,g)$ admits a global parallel null vector field, see \cite[Proposition 2]{BaumLarzLeistner14} for details.

If we now assume that $M$ is compact, the time-orientable cover is  compact and admits a global parallel null vector field, and hence, by the result in \cite{MehidiZeghib22}
is geodesically complete. This implies that $(M,g)$ is geodesically complete. \hfill $\qed$

\bibliographystyle{abbrv}
\bibliography{GEOBIB}
\end{document}